\documentclass[towside,a4paper,french,12pt]{article}
\usepackage[english]{babel}
\usepackage[T1]{fontenc}
\usepackage{babel}
\usepackage{psfrag}
\usepackage{amsmath}
\usepackage{amsmath,amsthm} 
\usepackage{lipsum}
\usepackage{amsfonts}
\usepackage{amssymb}
\usepackage{eufrak}
\usepackage{dsfont}
\usepackage{makeidx}
\usepackage{graphicx}
\usepackage{vmargin}
\usepackage{fancyhdr}
\usepackage{enumerate}
\usepackage{amsmath,amscd,latexsym}
\usepackage{amsfonts}
\usepackage{pst-node}
\usepackage{pst-plot}
\setmarginsrb{3cm}{3.3cm}{3.5cm}{4cm}
   {1.1cm}{1cm}{1cm}{1cm}

 \fancypagestyle{plain}{
    \fancyhead[R]{\thepage}
    \fancyhead[L]{}
    
    }
\newtheorem{theo}{Theorem}[section]
\newtheorem{defi}[theo] {Definition}%[section]
\newtheorem{prop}[theo]{Proposition}%[section]
\newtheorem{rem}[theo]{Remark}%[section]
\newtheorem{lem}[theo]{Lemma}%[section]
\newtheorem{coro}[theo]{Corollary}%[section]

\newenvironment{dem}{{\textbf{Proof}.}}{\hspace{\stretch{1}}$\square$\vspace{0,2 cm}}

\newcommand{\flech}[1]{\ensuremath{\searrow^{\hspace{-0.3 cm}{\small #1}}}}
\newcommand{\Aut}{{\rm Aut}}

\def\phi{\varphi}

\title{\textbf{Cylinders, multi-cylinders and the induced action of $\Aut(F_n)$}}
\author{Fedaa Ibrahim
\\
\\
\footnotesize{
Universit\'e d'Aix-Marseille,
Laboratoire d'Analyse, Topologie, Probabilit\'es}\\
\footnotesize{Avenue Escadrille Normandie-Niemen, 13397 Marseille, France}}

 \date{}

\begin{document}

\maketitle
\begin{abstract}
A cylinder  $C^1_u$  is the set of infinite words with fixed prefix  $u$.  A double-cylinder $C^2_{[1,u]}$ is ``the same'' for bi-infinite words.  We show that for every word  $u$  and any automorphism $\varphi$ of the free group $F$   the image $\varphi(C^1_u)$  is a finite union of cylinders.  The analogous statement is true for double cylinders.  We give  (a)  an algorithm, and  (b)  a precise formula which allows one to determine this finite union of cylinders.
\end{abstract}

\section{Introduction}

This paper goes back to a remark of a rather well known member of the ``Outer space" community, who some years ago during a talk in Bonn explained that rational currents are dense in the space of currents, but that, other than using this fact and a bit of approximation, she didn't know how to compute 
the image of a current under the induced action of an automorphisms $\varphi$ of a finitely generated free group $F$. 

By definition, a current $\mu$ is a measure on the {\em double boundary} $\partial^2 F$, i.e. the space $\partial F \times \partial F$ minus the diagonal.  The image measure $\varphi_*(\mu)$, of course, is simply the measure $\mu$ evaluated on the preimages of subsets of $\partial^2 F$ under the homeomorphisms induced by $\varphi$.  The problem, it turns out, is that even for the simplest sets in $\partial^2 F$, the so called {\em double cylinders} $C^2_{[u,v]}$
(see Definition \ref{5.3}), given by two distinct elements $u, v\in F$ and the choice of a basis $A$ of $F$, it is not at all evident how to describe $\varphi(C^2_{[u,v]})$ (or $\varphi^{-1}(C^2_{[u,v]})$).
For example, using the results of this paper, it is easy to give examples of double cylinders with $\varphi(C^2_{[u,v]})\neq C^2_{[\varphi(u),\varphi(v)]}$. 

Indeed, we prove here (see \S 5):

\begin{theo}
\label{main-thm}
Let $\varphi$ be an automorphism of the free group $F$ with finite basis $A$.
For any $u,v\in F$ with $u \neq v$ there exist finite sets $U,V \subset F$ such that:
$$\varphi (C^{2}_{[u,v]})=\dot{\bigcup \limits_{\substack {u_{i} \in U \\ v_{j}\in V}}}  \ C^{2}_{[u_{i},v_{j}]}$$
The sets $U$ and $V$ can be algorithmically derived from $u, v \in F$ and from the elements of $\varphi(A)$ and of $\varphi^{-1}(A)$, all expressed as reduced words in $A \cup A^{-1}$.
\end{theo}

To simplify the arguments, one considers first one-sided cylinders $C^1_w \subset \partial F$:
they too depend on the chosen basis $A$ of $F$, since one has to pass from the element $w \in F$ to the corresponding element of $F(A)$, by which we denote the set of reduced words in $A \cup A^{-1}$. One thus obtains $C^1_w$
as the set of all 
elements of $\partial F$ that are represented by
one-sided infinite reduced words in $A \cup A^{-1}$ which have $w$ as prefix. We also need to consider multi-cylinders $C^1_U = \bigcup \limits_{u \in U} C^1_u$ for finite subsets $U\subset F$. In \S 4 below we show:

\begin{theo}
\label{main-thm-1}
Let $\varphi$ be an automorphism of the free group $F$ with finite basis $A$.

\smallskip
\noindent 
(a) 
For any $u \in F(A)$ there exists a finite set $U \subset F(A)$ such that:
$$\varphi (C^1_u)=C^{1}_U$$

\smallskip
\noindent 
(b) A set $U$ as in statement (a) can be algorithmically derived from $u \in F(A)$ and from the words in the finite subsets $\varphi(A)$ and $\varphi^{-1}(A)$ of $F(A)$. Indeed, the equality in (a) is true for
$$U = \{\varphi(u')|_{S(\varphi)^2} \mid u' \in u|^k\}\, ,$$
with 
$k=S(\varphi)^4+S(\varphi)^3+S(\varphi)^2$, where $S(\varphi)$ is the maximal length of any $\varphi(a_i)$ or $\varphi^{-1}(a_i)$ among all $a_i \in A$, see \S 2.

Here
for any reduced word $w \in F(A)$ and any integer $l \geq 0$ we denote by $w|_l$ the word obtained from $w$ by erasing the last $l$ letters, and by $w|^l$ the set of reduced words obtained from $w$ by adding $l$ letters from $A \cup A^{-1}$ at the end of $w$.
\end{theo}

The set $U$ from the above Theorem \ref{main-thm-1} is not uniquely determined by $u$, $A$ and $\varphi$:  The set $U$ exhibited in part (b) is only one of infinitely many finite subsets $U' \subset F(A)$ which all satisfy the equality $\varphi (C^1_u)=C^{1}_{U'}$ from part (a).  

\smallskip

This non-uniqueness can be easily understood by considering the following two
typical examples, given by the pairs 
$U_1=\{ab, aba^{-1} \} , U_2=\{ab \} $ and by
$U_3=\{ aba^{-1},aba,abb \} , U_4=\{ab \} $, which satisfy $C^1_{U_1} = C^1_{U_2}$ and $C^1_{U_3} = C^1_{U_4}$.
The resulting ambiguity is resolved by the following proposition, which is proved below in \S 3:

\begin{prop}
\label{thm1-}
For every multi-cylinder $C^1_U$, determined by a finite set $U \subset F(A)$, there is a unique finite subset $U_{\min} \subset F(A)$ of minimal cardinality which determines the same multi-cylinder: 
$$C^1_{U_{\min}} = C^1_U$$ 
The set $U_{\min}$ can be derived algorithmically from $U$ by a finite sequence of elementary operations (of two types, illustrated by the two examples presented in the previous paragraph), each of which strictly decreases the cardinality.
\end{prop}

This enables us to define a map $\varphi^*_A$ on elements (and on finite subsets of $F(A)$) by associating to $u \in F(A)$ the minimal set $U_{\min}$ for the multi-cylinder $\varphi(C_u)$:  the set $U_{\min}$ can be derived algorithmically from any finite set $U \in F(A)$ as in Theorem \ref{main-thm-1}, with $C^1_U = \varphi(C_u)$.

\smallskip

We can thus reformulate and specify the main case of Theorem \ref{main-thm} slightly, by stating (see \S 5):

\begin{prop}
Let $u,v\in F(A)$ be such that none is prefix of the other.
Then one has  
$$\varphi (C^{2}_{[u,v]})= \dot{\bigcup \limits_{\substack {u_{i} \in \varphi^*_{A}(u) \\ v_{j}\in \varphi^*_{A}(v)}}}  \ C^{2}_{[u_{i},v_{j}]}$$
\end{prop}

The extra hypothesis in the last proposition is necessary since double cylinders behave properly under the action of $F$ on the indices (see Lemma \ref{5.6}), while for a single cylinder $C^1_u$ one has $w C^1_u = C^1_{w u}$ only if $u$ is not a prefix of $w^{-1}$.  For a general formula see Remark \ref{pas-encore-ecrit}.

\bigskip

\noindent
{\em Acknowledgments:} 

I would like to thank my two co-advisors Thierry Coulbois and Martin Lustig for constant encouragement and many helpful remarks 
throughout my time as graduate student; this paper is indeed part of my thesis \cite{FI}.  
My special thanks go also to Pierre Arnoux for having directed my attention towards the proof of Proposition \ref{f16}, and to Llu\'is Bacardit for some valuable literature remarks. Furthermore, I would like to thank Nicolas Bedaride for his efforts concerning the careful proofreading of \cite{FI}. Finally, the report of the referee was very helpful for producing the final version; I am also grateful for that.

\section{Notation, set-up and basic facts}
Throughout this paper we denote by $F$ a finitely generated non-abelian free group, and by $\varphi$ an automorphism of $F$.
We choose a basis $A$ of $F$ once and for all, which allows us to identify $F$ with the set $F(A)$ of finite reduced words in the elements of
$A$ and their inverses.
We denote by
$\partial F(A)$  the set of infinite reduced words:
$$\partial F(A)=\{x_{1} x_{2} x_{3} \cdots  \mid  x_{i} \in A\cup  A^{-1} , \ x_{i} \neq x_{i+1}^{-1} \}$$
The set $\partial F(A)$ is in a canonical bijective correspondence with the {\em end completion} 
$\partial  F$ of $F$. The latter also coincides with the {\em Gromov boundary} of $F$. The set $\partial F$ (and thus $\partial F(A))$ carries a topology; indeed it is homeomorphic to a Cantor set.  Every automorphism $\varphi$ of $F$ induces canonically a homeomorphism of $\partial F$, 
which for simplicity we denote also by $\varphi$.
For background and details about these classical facts see \cite{BH}.

The \emph{word length} of an element $w\in F(A)$ with respect to $A$ will be denoted by $|w|_A$ or simply by $|w|$.
We write $v\leq w$, if $v$ is a prefix (= initial subword) of $w$, and we write $v<w$ if in addition one has $|v| < |w|$.
This puts a partial ordering on $ F$ (which heavily depends on $A$).
The longest prefix common to elements $w_{1}$ and $w_{2}$ of $F(A) \cup \partial F(A)$ is denoted
$w_{1} \wedge w_{2}$. 
One has $|w_1^{-1} \wedge w_2|=0$ if and only if the product $w_1 w_2$ is reduced; in this case we denote  $w_1 w_2$ by $w_{1} \cdot w_{2}$.

The \emph{size} of an automorphism $\varphi \in \Aut (F)$ (with respect to $A$) is defined by
$$S( \varphi ) := S_A( \varphi ) := \max \limits_{a\in A\cup  A^{-1}} \{| \varphi (a) |, | \varphi^{-1} (a) | \}$$

We obtain directly from this definition:

\begin{lem}\label{f2,1}
For any $w \in F(A)$ and any $\varphi \in \Aut (F)$ one has:
$$\frac{|w|}{ S(\varphi)} \,\,\, \leq \,\,\, |\varphi (w)| \,\,\, \leq \,\,\, |w| \cdot S(\varphi)$$
\qed
\end{lem}

The following is a classical result of D. Cooper, see \cite{MR916179}. 

\begin{prop}\label{f2,3} 
Let $\varphi$ be an automorphism of the finitely generated free group $F$, and let $A$ be a basis of $F$.
Then there exists a constant $C\geq 0$ such that for any elements $u,v \in F$ one has:
$$0\leq |\varphi (u)|_{A} + |\varphi (v)|_{A}-|\varphi (uv)|_{A} \leq C$$
The smallest such constant $C$ will be denoted by $C(\varphi)$.
\end{prop}

In the literature the above proposition is sometimes referred to as ``bounded cancellation lemma''.
It follows directly from this proposition that the analogous statement, i.e. the upper bound on the possible cancellation, remains true if $u^{-1}$ or $v$ (or both) are replaced by elements from $\partial F$, i.e. by infinite words.

\begin{rem}
\rm
In \cite{MR916179} it has been shown that for any $\varphi \in \Aut(F)$ the constant $C(\varphi)$ is always bounded above by $S(\varphi)^2$.
\end{rem}

\begin{defi}
\rm
Let $w= a_{1}\cdots a_{r}\in F(A)$. For any integer $k \geq 0$ we define:

\smallskip
\noindent
(1)
 $w|_{k}=a_{1} \cdots a_{r-k}$ \,\, (if $k \leq r$), and
 
\smallskip
\noindent
(2)
$w|^{k}=\{v \mid w < v \,\, {\rm and} \,\,  |v| = |w|+k \}$
\end{defi}

From this definition we obtain directly, for any 
$u\in F(A)$ and any integers $m,n \geq 0$ with $k=m+n$,
that 
$u|^k= \bigcup \limits_{v\in u|^m} v|^n$.

%%%%%%%%%%%%%%%%%%%%%%%
\section{Cylinders and multi-cylinders}

It is crucial in this section that one distinguishes between elements of the free group $F$, with basis $A = \{a_1, \ldots, a_n\}$, and reduced words in the $a_i$ and $a_i^{-1}$ which are used to represent these elements. We denote the set of reduced words by $F(A)$.

Similarly, we denote by $\partial F(A)$ the set of infinite reduced words $X = x_1 x_2 \ldots$ in $A \cup A^{-1}$ which are used to represent the elements of the Gromov boundary $\partial F$.

We will denote in this section by $\mathds{U}$ the set of all  finite subsets of $F(A)$.

\begin{defi}
\label{3.0}
\rm
For any  $u\in F(A)$ we define
$C^1_u= \{X\in \partial{F}(A) \mid  u < X \}$. The set $C^1_u$ is called the $\emph{cylinder}$ defined by $u$ (and by $A$).
\end{defi}

\begin{rem}\label{3,1}
\rm
Let $u,v \in F(A)$.  Then from the definition of $C^1_u$ one derives directly$:$
\begin{enumerate}[(1)]
\item If $C^1_u=C^1_{v}$ then $u=v$.
\item If $C^1_u \cap C^1_{v}\not=\emptyset$ then one has $v\leq u$ and thus $C^1_u \subseteq C^1_v$,  or  else  $u\leq v$ and thus $C^1_v\subseteq C^1_u$.
\item For any integer $k\geq 0$ one has $C^1_u=\dot {\bigcup \limits_{u_i\in u|^k}} C^1_{u_i}$.
\end{enumerate}
\end{rem}

From parts (1) and (2) of Remark \ref{3,1} we obtain directly$:$

\begin{lem}\label{3,2}
Given $u,u'\in F(A)$ with  $|u|=|u'|$, then either $C^1_u \cap C^1_{u'}=\emptyset$,  or else $u=u'$ and thus $C^1_u=C^1_{u'}$. \qed
\end{lem}

\begin{defi}
\rm
For any subset $U \subset F(A)$ we will denote by $C^1_{U}\subset \partial F$ the union of all cylinders $ C^1_{u}$ with $u \in U$$:$
$$C^1_U= \bigcup \limits_{u_i \in U} C^1_{u_i}$$
\end{defi}

From Lemma \ref{3,2} we obtain directly$:$

\begin{lem}\label{3,4}
Let $k\in\mathds{N}$, $U\subset F(A)$ and $|u_i|=k$ for all $u_i\in U$. Then one obtains a disjoint union: 
$$C^1_U=\dot {\bigcup \limits_{u_i\in U}} C^1_{u_i}$$ 
\qed
\end{lem}

Recall that $\mathds{U}$ denotes the set of all finite subsets of $F(A)$.

\begin{lem}\label{f11,2}
Let $k\in\mathds{N}$ and $U,U'\in \mathds{U}$, and assume for all $u\in U\cup U'$ that $|u|=k$. Then we have $C^1_U=C^1_{U'}$ if and only if $U=U'$.
\end{lem}

\begin{dem}
If $U=U'$ then clearly one has $C^1_U=C^1_{U'}$.
Conversely, from the hypothesis  $|u|=k$ for all $u\in U\cup U'$ we obtain, by Lemma \ref {3,4}, that  $C^1_U=\dot {\bigcup \limits_{u_i\in U}} C^1_{u_i}$ and  $C^1_{U'}=\dot {\bigcup \limits_{u'_j\in U'}} C^1_{u'_j}$. Thus, if $C^1_U = C^1_{U'}$, we obtain  $\dot {\bigcup \limits_{u_i\in U}} C^1_{u_i} =\dot {\bigcup \limits_{u'_j\in U'}} C^1_{u'_j}$. From Lemma \ref{3,2} we deduce that for any $C^1_{u_i} \subset C^1_U$ there exists a unique $C^1_{u'_j}\subset C^1_{U'}$ with $C^1_{u_i}=C^1_{u'_j}$ and thus $u_i=u'_j$ (by Remark \ref{3,1} (1)). This shows  $U\subset U'$, and from the symmetry between $U$ and $U'$ we obtain $U=U'$.
\end{dem}

We define now an ``elementary'' relation $\searrow$ on $\mathds{U}$ as follows$:$\

\begin{defi}
\rm
For any $U_1,U_2 \in \mathds{U}$ we write $U_1\searrow U_2$
if one of the following conditions is satisfied$:$
\begin{enumerate}
  \item  There are distinct elements $u_i,u_j \in U_{1}$ with $u_{i} < u_{j}$ such that $U_{2}=U_{1}\smallsetminus \{u_{j}\}$. In this case we sometimes specify the notation $U_1\searrow U_2$ to  $U_2  \flech{(1)}  U_1$.
  \item There exists an element $u\in F(A) \smallsetminus U_1$ with   $u|^{1} \subset U_{1}$, and one has $U_{2}= (U_{1} \smallsetminus u|^{1}) \cup \{u\}$. In this case we
write sometimes $U_2 \flech{(2)} U_1$.
\end{enumerate}
\end{defi}

For example, let $F$ be a free group with base $A=\{a,b\}$, and let $U= \{aba, abab, bba, bbb, bba^{-1}\}$. Then for $U_1=\{aba, bba, bbb,bba^{-1}\}$ we have $U \flech{(1)} U_1$, and for $U_2=\{aba,bb\}$ we obtain $U_1 \flech{(2)} U_2$.

\begin{rem}\label{3,5}
\rm
It is clear that the relation $\searrow$  strictly decreases  the cardinality of the given set $U$:
$$U \searrow U' \,\,\, \Longrightarrow \,\,\, \#U > \#U'$$
\end{rem}

\begin{defi}
\label{relation}
\rm
For any $U, U' \in \mathds{U}$ we
write 
$U\sim U'$
if there exists a finite sequence $U_1=U, U_2,\cdots ,U_n =U'$ of elements of $\mathds{U}$, with $U_{i} \searrow  U_{i+1}$ or $U_{i+1}\searrow U_{i}$ for all $1 \leq i \leq n-1$.

In other words : The relation $\sim$ is the equivalence relation on $\mathds{U}$ generated by the elementary relation
$\searrow$ .
\end{defi}

\begin{defi}
\rm
We say that $U\in \mathds{U}$ is reduced if and only if there is no
 $U'\in \mathds{U}$ with $U\searrow U'$.
\end{defi}

\begin{rem}\label{f11,0}
\rm
(a)
For any $U\in \mathds{U} $ there exists a reduced set $U' \in \mathds{U}$ with $U\searrow \dots \searrow U'$.
{This follows directly from the finiteness of $U$ and from Remark \ref{3,5}}.

\smallskip
\noindent
(b)  However, it is a priori not clear that the reduced set $U'$ depends only on $U$ and not on the particular way how one choses the reduction $U\searrow \dots \searrow U'$. To show that in each equivalence class $[U]_\sim$ there is precisely one reduced set $U'$ is the goal of the rest of this section.
\end{rem}

\begin{lem} \label{f11,1}
(a)  Let $U,U' \in \mathds{U}$ and assume $U\searrow U'$. Then we have $C^1_U=C^1_{U'}$.\\
(b)  In particular,  if $U\sim U'$ then one has $C^1_U=C^1_{U'}$.
\end{lem}

\begin{dem}
(a)  From the above definition of $\searrow$ we distinguish two cases:

\noindent
(1) If $U\flech{(1)} U'$ then there exist $u_1,u_2\in U$ with $u_1< u_2$ and $U'=U\smallsetminus \{u_2\}$. Thus one has $U=U'\cup \{u_2\}$, and thus $C^1_U=C^1_{U'}\cup C^1_{u_2}$. But $C^1_{u_2} \subset C^1_{u_1} \subset C^1_{U'}$, so that $C^1_U=C^1_{U'}$.

\smallskip
\noindent
(2)
If $U\flech{(2)} U'$ then there exists $u\in F(A)$ with $u\not \in U$, $u|^1 \subset U$
and $U'= (U\smallsetminus u|^1 )\cup \{u\}$. Thus we have $C^1_U=C^1_{U'\smallsetminus \{u\}} \cup C^1_{u|^1}$ and $C^1_{U'}=C^1_{U\smallsetminus u|^1} \cup C^1_{u}$. From Remark \ref{3,1} (3) one has $C^1_u=C^1_{u|^1}$, so that the last two equalities give $C^1_U \supset C^1_{U'}$ and $C^1_{U'} \supset C^1_{U}$, and thus $C^1_U=C^1_{U'}$.

\smallskip
\noindent
(b) This is a direct consequence of (a), by the definition of $\sim$.
\end{dem}

We now define another elementary relation $\nearrow $ which allows us to extend a set $U_1$ to a larger set $U_2$$:$ 

For any $U_1, U_2 \in \mathds U$ we write
$U_1\nearrow U_2$ if 
{$u\in U_1 $ and  $U_2= U_1 \cup u|^1 \smallsetminus \{u\}$.}

\begin{rem}\label{f11,2}
\rm
(a) We observe that  $U_1\nearrow U_2$ does not necessarily imply that $U_2 \flech{(2)} U_1$. For example, if $U_1=\{b,ba\}$ and $U_2=\{ba, bb, ba^{-1}\}= b|^1$ then we have $U_1\nearrow U_2$ and $U_2 \flech{(2)}  \{b\} 
\subsetneqq U_1$.

\smallskip
\noindent
(b) If $U_1 \nearrow U_2$ then one has $U_1\sim U_2$. To see this, we observe from $U_1 \nearrow U_2$ that  there exists $u\in U_1$ such that  $u\not\in U_2$ and $ u|^1\subset U_2$. 
Now we apply  \flech{(2)} to obtain $U_2 \flech{(2)} U_2'$, where $U_2' =\{U_2 -u|^1\}\cup \{u\}$. Thus all elements of $U_2-U'_2$ are contained in $u|^1$. Since $u\in U'_2$, a multiple application of \flech{(1)} yields $U_2 \flech{(1)}  \cdots \flech {(1)} U'_2$. This implies  $U_1\sim U'_2$.

\smallskip
\noindent
(c)  In particular, by Lemma \ref{f11,1} (b), if $U_1\nearrow U_2$ then $C^1_{U_1}= C^1_{U_2}$.
\end{rem}
\begin{prop}\label{f12}
For all $U,U' \in \mathds{U}$  one has$:$ $$C^1_{U}=C^1_{U'} \,\,\, \Longleftrightarrow \,\,\, U\sim  U'$$
\end{prop}

\begin{dem}
If $U\sim U'$ then by  Lemma \ref{f11,1} (b) we have $C^1_U=C^1_{U'}$.
For the converse direction assume $C^1_U=C^1_{U'}$. Let $k=\max\{|u|\mid u\in U\cup U'\}$. We set $U_0=U$ and define iteratively $U_{i+1} $ from $U_i$ by postulating
$$ U_{i+1}=
(U_i  \smallsetminus \{u\})\cup u|^1 
$$ 
for some $u\in U_i $ with $ |u| <k $.
Then one obtains $U=U_i\nearrow  U_{i+1} \nearrow  U_{i+2} \nearrow \cdots \nearrow U_n$, where for all $ v\in U_n$ one can assume $|v|=k$. By part (b) of Remark \ref{f11,2} we obtain $U\sim U_n$ and thus  $C^1_U=C^1_{U_n}$.

We do the same for $U'$ to find $U'=U'_0\nearrow U'_{1} \nearrow \cdots \nearrow U'_m$, where for all $ v'\in U'_m$ one has $|v'|=k$. Again we obtain  $U'\sim U'_m$ and thus $C^1_{U'}=C^1_{U'_m}$. But we assumed $ C^1_U=C^1_{U'}$, which gives $C^1_{U_n}=C^1_{U'_m}$ and thus, by Lemma \ref{f11,2} , $U_n=U'_m$. This gives  $U\sim U_n=U'_m \sim U'$ and hence $U\sim U'$.
\end{dem}

\begin{defi}\label{12}
\rm
For any subset $B \subset \partial F(A)$ we define
$$U^*(B)=\{ u\in F(A) \mid C^1_u\subset B \,\,{\rm and}\,\, C^1_{u|_1}\not \subset B \}\, .$$
For $U\in \mathds{U}$ we write $U^* := U^*(C^1_U)\in \mathds{U}$.
\end{defi}

\begin{rem}\label{f12,000}
\rm
From Definition \ref{12} we obtain directly: \\
(a)   If $U,V\in\mathds{U}$, with $C^1_U=C^1_V$, then $U^*=V^*$.\\
(b)   For all $U\in \mathds{U}$ we have $C^1_{U^*} \subset C^1_U$.\\
(c)   For all $U\in \mathds{U}$ one has $(U^*)^* = U^*$.
\end{rem}

\begin{lem}\label{f12,0}
For any $U\in \mathds{U}$ one has $C^1_{U^*}=\dot{\bigcup \limits_{u\in U^*}}  C^1_u$.
\end{lem}

\begin{dem}
If, by way of contradiction, we assume  $C^1_{U^*}\not=\dot{\bigcup \limits_{u\in U^*}}  C^1_u$, then there exist $u_1,u_2\in U^*$, $u_1\not= u_2$, with $C^1_{u_1}\cap C^1_{u_2}\not= \emptyset$. \   By part (2) of Remark \ref{3,1}   one has  $u_1<u_2\,$ or $\,u_2<u_1$ and thus $u_1 \leq u_2|_1$ or $u_2\leq u_1|_1$. This implies  $C^1_{u_2|_1}\subset C^1_{u_1} \subset C^1_U$ or $C^1_{u_1|_1}\subset C^1_{u_2} \subset C^1_U$, which contradicts the assumption  $u_1,u_2\in U^*$. Hence we have proved $ C^1_{U^*}=\dot{\bigcup \limits_{u\in U^*}}  C^1_u$.
\end{dem}

\begin{lem}\label{f12,00}
For each $U\in \mathds{U} $ there is no $U'\sim U$ with $U'\subsetneqq  U^*$.
\end{lem}

\begin{dem}
From Lemma \ref{f12,0} we know $C^1_{U^*}=\dot{\bigcup \limits_{u\in U^*}}  C^1_u$, and from Remark \ref{f12,000} (b) we have $C^1_{U*} \subset C^1_U$. On the other hand, $U'\sim U$ implies by Proposition \ref{f12} the equality $C^1_U=C^1_{U'}$ and thus $C^1_{U^*} \subset C^1_{U'}$. As a consequence, one deduces from $U'\subset U^*$ that $\dot{\bigcup \limits_{u\in U^*}}  C^1_u =\dot{\bigcup \limits_{u\in U'}}  C^1_u$, which implies $U'=U^*$, since every $C^1_u$ is non-empty.
\end{dem}

\begin{lem}\label{f12,1}
If $U\in \mathds{U} $ is reduced, then one has $U= U^*$.
\end{lem}

\begin{dem}
By way of contraction assume $U\not =U^*$. By Lemma \ref{f12,00} this implies that $U-U^*$ is non-empty. Let  $n=\max\{|u|\mid u\in U-U^*\}$, and let  $u\in U-U^*$ with $|u|=n$. 
By definition of $U^*$ we have that $C^1_{u\mid_1}\subset C^1_U$, so that one of the following three properties must hold:

(1)   $u|_k \in U$ for some $k \geq 1$.

(2)   $u|_1\big |^1 \subset U$.

(3)   $u|_k \not\in U$ for all $k \geq 1$, and there
exists $v \in  u|_1\big |^1$ (i.e.  $|v|=n$) with $v\not \in U$.\\ 
The cases (1) and (2) are impossible because $U$ is reduced and $u\in U$.
 In case (3), since  $C^1_v\subset C^1_{u|_1} \subset C^1_U$, there exists $v' \in  v\big |^k$, with $k\geq1$,   $|v'|=n+k$, $v' \in U$ and $C^1_{v'}\subset C^1_U$. We deduce  $C^1_{v'|_1}\subset C^1_{u|_1} \subset C^1_U$, and thus $v'\in U-U^*$$:$  This contradicts the definition of $n$ above because $|v'|>n$.
 \end{dem}

\begin{prop}\label{f12,2}
(a) For every $U\in \mathds{U} $ there is precisely one reduced set $U_{\min} \in \mathds U$ with $U_{\min} \sim U$.  

\smallskip
\noindent
(b) In particular, one has $U_{\min} = U^*$ and $C^1_U=C^1_{U_{\min}}=C^1_{U^*}$,
and this is the disjoint union of all $C^1_u$ with $u \in U_{\min}$.
 \end{prop}
 
 \begin{dem}
Let $U' \in \mathds U$ be a reduced set with $U\sim U'$.  By Remark \ref{f11,0} (a) such a set $U'$ exists.
By Proposition \ref{f12} we have $C^1_U=C^1_{U'}$ and thus $U^* = U'^*$. As $U'$ is reduced, by Lemma \ref{f12,1} we have $U'=U'^*$ and thus $U' = U^*$. This shows the uniqueness of the set $U' =: U_{\min}$, as well as the equalities stated in claim (b).
\end{dem} \\

We now obtain Proposition \ref{thm1-} stated in the Introduction as an immediate consequence of Remark \ref{f11,0}, Proposition \ref{f12} and Proposition \ref{f12,2}.

%%%%%%%%%%%%%%%%%%%%%%%%%%%%
\section{The $\varphi$-image of a cylinder $C^{1}_{w}$}
The objective of this section is to determine the image of any cylinder $C^1_w$, with 
$w \in F(A)$, under a given automorphism $\varphi$ of the free group $F(A)$. We will see that there exists a finite set $U\subset F(A)$  of words in $A$  such that 
$$\varphi (C^1_w)=\dot {\bigcup \limits_{u\in U}} C^1_u$$  
In this section we will first prove  the existence of such a finite set $U$, and in a second step  we will define  an algorithm  that determines $U$, for any given word $w\in F(A)$ and any automorphism $\varphi $ of $F(A)$ (given by the finite set of words $\varphi(a_i)$ for any $a_i \in A$).

\begin{rem}\label{f15}
\rm
Given $w \in F(A)$, we first note that in general one has$:$
$$\varphi (C^{1}_{w}) \neq C^{1}_{\varphi(w)}$$
For example, let $F(a,b)$ be the free group with base $\{a,b\}$,
and let $\varphi\in \Aut (F(a,b))$, given by$:$
$$a\mapsto aba \ \ , \ \ \ b\mapsto ba $$
We consider  $w =ba $ and obtain $\varphi(w)=baaba$, as well as \\
\centerline{$C^{1}_{w}=\{baz_{1}z_{2} \cdots \mid z_{1}\in \{a,b, b^{-1}\}, z_i \in \{a,b,a^{-1}, b^{-1}\} \smallsetminus \{z^{-1}_{i-1}\} \ \forall i \geq 2  \}$} 
and\\
\centerline{
$C^{1}_{\varphi(w)}=\{baabaz_{1}z_{2} \cdots \mid z_1\in \{a,b, b^{-1}\}, z_{i}\in \{a,b,b^{-1},a^{-1}\}\smallsetminus \{z_{i-1}^{-1}\} \ \forall i\geq 2  \}$.}
Then for
$W=bab^{-1}a^{-1}a^{-1}a^{-1}a^{-1} \cdots \in C^1_w$ we obtain\\
\centerline{$\varphi (W)=bab^{-1} a^{-1}  a^{-1} b^{-1}a^{-1} a^{-1} b^{-1}a^{-1} a^{-1} b^{-1}a^{-1} \cdots \in \varphi(C^1_w)$\ ,}
and we observe $\varphi(W) \not \in C^1_{\varphi(w)}$, which implies $\varphi (C^{1}_{w}) \neq C^{1}_{\varphi(w)}$.
\end{rem}

We'd like to thank P.~Arnoux for having pointed out to us that a proof of the following statement should be possible along the lines given below in the proof.

\begin{prop}\label{f16}
For any $\varphi \in \Aut (F)$ and   \ $w \in F(A)$ there is a finite set  
$U \subset F(A)$ 
such that 
$$\varphi (C_{w}^{1})= \bigcup \limits_{u_{i} \in U} \ C^{1}_{u_{i}} $$
\end{prop}
\begin{dem}
With respect to its natural topology (see \S 2) the space $\partial F$ is compact, and for any $u \in F(A)$ the cylinder $C_u^1$ is open and compact.
 Since every $\varphi \in \Aut (F) $ induces a homeomorphism on $\partial F$, for any $u\in F(A)$ the image set $\varphi (C_u^1) $ must also be open and compact. 
Thus, since the set $\{C^1_u \mid u \in F \}$ constitutes a basis of the topology of $\partial F$, it follows from $\varphi (C_u^1) $ open that there is a (potentially infinite) family of $C^1_{u_i} \subset \varphi (C_u^1) $ which covers all of $\varphi (C_u^1) $.  By the compactness of the latter we can extract a finite subfamily $\{C^1_{u_i} \mid u \in U \}$ which still covers $\varphi (C_u^1) $, while each $C^1_{u_i}$ remains a subset of $\varphi (C_u^1) $. This proves the claim.
\end{dem}

It should be noted that the above proof of Proposition \ref{f16} has no algorithmic value. Indeed, it does not even allow us to find $U$ by trial and error 
(unless one first derives an algorithm that verifies the equality of Proposition \ref{f16} for any given $\varphi, w$ and $U$).

\begin{lem} \label{f17}
Let $ \varphi \in \Aut(F)$ \ and \ $ w \in F(A)$ with $|w|\geq  S(\varphi)\cdot C(\varphi)$.\\ Then one has$:$
$$\varphi (C^{1}_{w})\subset C^{1}_{\varphi (w)|_{C(\varphi)}}$$
\end{lem} 

\begin{dem} For all $Z\in C_{w}^{1}$ there exists $X\in \partial F(A)$ such that  $Z= w \cdot X$ and hence 
$\varphi(Z) \in \varphi( C_{w}^{1})$ and \  $\varphi(Z)= \varphi(w) \varphi(X)$.
By the definition of $S(\varphi )$ (see \S 2) we have 
$|\varphi(w)|\geq \frac{|w|}{ S(\varphi)}$, 
and by assumption we know $|w|\geq S(\varphi) \cdot C(\varphi)$,
so that $|\varphi(w)|\geq C(\varphi)$. Thus we can
decompose $\varphi(w) = w_{1}\cdot w_{2}  $, where $|w_{2}|=C(\varphi )$ and $w_1=\varphi(w)|_{C(\varphi)}$.
The cancelation between $\varphi(w)$ and $\varphi(X)$ is bounded by $C(\varphi)$ 
(see Proposition \ref{f2,3} and the subsequent paragraph),
so that 
for some decomposition $w_2 = w'_2 \cdot w''_2$ 
we obtain  $\varphi(Z)=w' \cdot X'$ with $w'=w_1\cdot w'_2$ and
$\varphi(X)=w'^{-1}_{2} \cdot X'$.
This shows $\varphi(Z)\in C^1_{w'} \subset C^{1}_{w_{1}}$, which in turn proves  $\varphi(C^{1}_{w}) \subset C^{1}_{\varphi(w)|_{C(\varphi)}}$.
\end{dem}
\begin{prop} \label{f18}
Let $u,u' \in F(A)$, and assume: 
\begin{enumerate}
                            \item $u \leq u'|_k $ for $k=S(\varphi)\cdot C(\varphi)+C(\varphi^{-1})$
                            \item $|\varphi(u')|\geq S(\varphi)\cdot C(\varphi^{-1})+C(\varphi)$
                          \end{enumerate}
Then one has: $$C^1_{\varphi (u')|_{C(\varphi)}} \subset \varphi(C^1_u)$$
\end{prop}
\begin{dem}
From hypothesis 1. we obtain that $|u'|\geq  S(\varphi)\cdot C(\varphi)$, and thus we deduce from  Lemma \ref{f17} that \\ (I) $\varphi (C^1_{u'}) \subset C^1_{\varphi(u')|_{C(\varphi )}}$.\\
As a direct consequence we obtain that \\ 
(II)  $C^1_{u'}=\varphi^{-1} (\varphi (C^1_{u'}) ) \subset \varphi^{-1} (C^1_{\varphi(u')|_{C(\varphi)}} )$.\\
Now we apply hypothesis 2. to obtain $ | \varphi(u')|_{C(\varphi)} | \geq S(\varphi)\cdot C(\varphi^{-1})$.
This allows us to again apply Lemma \ref{f17}, with $w=\varphi(u')|_{C(\varphi)} $
and with $\varphi^{-1}$ instead of $\varphi$, to obtain \\
(III)  $\varphi^{-1}(C^1_{\varphi(u')|_{C(\varphi) }}) \subset C^1_{\varphi^{-1}(\varphi(u')|_{C(\varphi)})|_{C(\varphi^{-1})}}$.\\
From (II) and (III) we deduce \\
(IV)  $C^1_{u'}\subset C^1_{\varphi^{-1}(\varphi(u')|_{C(\varphi)})|_{C(\varphi^{-1})}}$, \\ which is equivalent to\\
(V)  $\varphi^{-1}(\varphi(u')|_{C(\varphi)})|_{C(\varphi^{-1})}\leq u'$.\\
By hypothesis 2. we can write $\varphi(u'):=u''\cdot u'''$ with $|u'''|=C(\varphi)$ and  $ u''=\varphi(u')|_{C(\varphi)}$. We calculate
\begin{equation*}
\begin{array}{lll}
\displaystyle\vert u'\vert&=&\vert \varphi^{-1}(u''\cdot u''')\vert\\
\displaystyle &\leq& \vert \varphi^{-1}(u'')\vert +\vert \varphi^{-1}(u''')\vert\\
\displaystyle &\leq & |\varphi^{-1}(u'')|+ S(\varphi)\cdot C(\varphi)
\end{array}
\end{equation*}
and thus obtain  $$|\varphi^{-1}(u'')|-C(\varphi^{-1} )\geq |u'|-S(\varphi)\cdot C(\varphi)- C(\varphi^{-1})\, .$$
As $u''=\varphi(u')|_{C(\varphi)}$, we can rewrite the last inequality as$:$
$$|\varphi^{-1} \left( \varphi(u')|_{C(\varphi)}\right)|-C(\varphi^{-1})\geq |u'|-S(\varphi)\cdot C(\varphi)-C(\varphi^{-1})$$
But $$\big| \varphi^{-1}\left( \varphi(u')|_{C(\varphi)}\right)|_{C(\varphi^{-1})}\big|=\big| \varphi^{-1}\left( \varphi(u')|_{C(\varphi)}\right)\big|-C(\varphi^{-1})$$
so that we obtain $\big| \varphi^{-1}\left( \varphi(u')|_{C(\varphi)}\right)|_{C(\varphi^{-1})}\big|\geq |u'|-k$. Hence we obtain from (V) that 
$u'|_k \leq \varphi^{-1}(\varphi(u')|_{C(\varphi)})|_{C(\varphi^{-1})} $, and thus from hypothesis 1. that 
$$u\leq\varphi^{-1}(\varphi(u')|_{C(\varphi)})|_{C(\varphi^{-1})}.$$
This is equivalent to $C^1_{\varphi^{-1}\left(\varphi(u')|_{C(\varphi)}\right)|_{C(\varphi^{-1})}} \subset C^1_u$.
From (III) we then deduce  that $\varphi^{-1}(C^1_{\varphi(u')|_{C(\varphi)}})\subset C^1_u$, which is equivalent to $$C^1_{\varphi(u')|_{C(\varphi)}} \subset \varphi(C^1_u)$$
\end{dem}

\begin{prop}\label{f19}
Let $u \in F(A)$ with  $|u|\geq  S^2(\varphi) C(\varphi^{-1})-C(\varphi^{-1})$, and let $k= S(\varphi)\cdot C(\varphi)+C(\varphi^{-1})$.
Then one has: $$\varphi (C^1_u)=\bigcup \limits_{u' \in u|^k} \ C^{1}_{\varphi(u')|_{C(\varphi)}}$$
\end{prop}
\begin{dem}
For all $u'\in u|^k$ one has $|u'|\geq k \geq S(\varphi) \cdot C(\varphi)$. Thus by Lemma \ref{f17}
we obtain $\varphi(C^1_{u'}) \subset C^1_{\varphi(u')|_{C(\varphi)}}$. Recall from 
part (3) of Lemma \ref{3,1}
that  $C^1_u = \bigcup \limits_{u' \in u|^k} C^{1}_{u'}$, which gives  $\varphi(C^1_u)=\varphi({\bigcup \limits_{u' \in
u|^k}} C^{1}_{u'})={\bigcup \limits_{u' \in u|^k}} \varphi
(C^{1}_{u'})$, so that one obtains \\
1.  \ \ \   $\varphi(C^1_u)\subset \bigcup \limits_{u' \in u|^k} C^1_{\varphi(u')|_{C(\varphi)}}$. \\

On the other hand, the hypothesis  $|u|\geq S^2(\varphi)C(\varphi^{-1} )-C(\varphi^{-1})$ is equivalent to
  $$|u|\geq S(\varphi)\left( S(\varphi)C(\varphi^{-1})+C(\varphi) \right) -S(\varphi)C(\varphi)-C(\varphi^{-1}),$$
which gives by $|u'|=|u|+k$ the inequality  \\
$|u'|\geq S(\varphi)\left( S(\varphi)C(\varphi^{-1})+C(\varphi) \right) -S(\varphi)C(\varphi)-C(\varphi^{-1})+S(\varphi)C(\varphi)+C(\varphi^{-1})$\\
$= S(\varphi)\left( S(\varphi)C(\varphi^{-1})+C(\varphi) \right)$.\\
Since $|\varphi(u')|\geq \frac{|u'|}{S(\varphi)}$ we obtain  $|\varphi(u')|\geq S(\varphi)C(\varphi^{-1})+C(\varphi)$.\\
Thus we can now apply Proposition \ref{f18}, to obtain 
$C^{1}_{\varphi(u')|_{C(\varphi)}} \subset \varphi (C^{1}_u)$ for all $u'\in u|^k$, so that one has \\
2. \ \ \   $\bigcup \limits_{u' \in u|^k} C^1_{\varphi(u')|_{C(\varphi)}} \subset
\varphi(C^1_u)$.\\
From 1. and 2. together we derive
$$\varphi(C^1_u)=\bigcup \limits_{u' \in u|^k} C^1_{\varphi(u')|_{C(\varphi)}}$$
\end{dem}

\begin{coro}
\label{4.7}
Let $k=k_1+k_2$, with $k_1=S^2(\varphi)C(\varphi^{-1})-C(\varphi^{-1})$ and \\ $k_2=S (\varphi)C(\varphi)+C(\varphi)$. Then for all $u\in F(A)$ we have
$$\varphi (C^1_u)=\bigcup \limits_{u' \in u|^k} \ C^{1}_{\varphi(u')|_{C(\varphi)}}$$
\end{coro}

\begin{dem}
For any  $v\in u|^{k_1}$ we have   $|v|\geq S^2(\varphi)C(\varphi^{-1})-C(\varphi^{-1})$.
Thus we can apply Proposition \ref{f19} to get
\begin{equation}
 \varphi (C^1_v)=\bigcup \limits_{u' \in v|^{k_2}} \ C^{1}_{\varphi(u')|_{C(\varphi)}}  
 \end{equation}
 Recall from part (3) of Remark \ref{3,1} that $C^1_u= \bigcup \limits_{v \in u|^{k_1}} \ C^{1}_{v} $ and thus $\varphi(C^1_u)= \bigcup \limits_{v \in u|^{k_1}} \ \varphi (C^{1}_{v}) $, so that we can deduce from equality (1):
 $$\varphi(C^1_u)=\bigcup \limits_{v \in u|^{k_1}} \big( \bigcup \limits_{u' \in v|^{k_2}} \ C^{1}_{\varphi(u')|_{C(\varphi)}} \big) $$
 Since $u|^k=u|^{k_1+k_2}$ 
 this is equivalent to
 $$\varphi (C^1_u)=\bigcup \limits_{u' \in u|^k} \ C^{1}_{\varphi(u')|_{C(\varphi)}}$$
\end{dem}

\begin{rem}
\label{alternative-approaches}
\rm
There are several alternative approaches to determine the image of a cylinder $C^1_u$ under an automorphism $\varphi$.  We briefly describe here two of them:

\medskip
\noindent
(a)  
Since every automorphism $\varphi$ of $F$ is a product of elementary automorphisms, one obtains a proof by induction over the length of such a product if one shows that for every elementary automorphism the image of a cylinder is a finite union of cylinders, and that those can be computed algorithmically. For permutations or inversions of the generators this is trivial; for elementary Nielsen automorphisms one has to work a little bit, but it is still not very difficult. On the other hand, this method doesn't permit one to describe $\varphi(C^1_u)$ by a closed formula as given in Corollary \ref{4.7}.

\medskip
\noindent
(b)  
Passing from $u \in F(A)$ to $u|^k$ for large $k$ is computationally rather an effort, so that the formula exhibited in Corollary \ref{4.7} is perhaps sometimes not very practical. We will thus sketch now a variation of the same basic approach, which has the advantage of being computationally more efficient (and also avoids some of the lengthly computations from above, after Lemma \ref{f17}):

\noindent
1. In a first step we pass from $u$ to some $u|^k$, but we pick the smallest possible $k \geq 0$ such that any $w \in u|^k$ satisfies the hypothesis of Lemma \ref{f17}.  This gives us a finite collection $W$ of words $w_i$ such that 
$\varphi(C_{u}^{1}) \subset \bigcup \limits_{w_i \in W} C_{w_{i}}^{1}$.

\noindent
2. We now prolong again every $w_i \in W$ to some $w_i|^{k_i}$, where $k_i \geq 0$ is chosen minimally to achieve two goals:

(i)  We can again apply Lemma \ref{f17} to any $u_j \in w_i|^{k_i}$, but this time with $\varphi^{-1}$ instead of $\varphi$. This gives $\varphi^{-1}(C^1_{u_j}) \subset C^1_{\varphi^{-1}(u_j)|_{C(\varphi^{-1})}}$.

(ii)  For any $u_j \in w_i|^{k_i}$ the word $\varphi^{-1}(u_j)|_{C(\varphi^{-1})}$ is not a prefix of $u$.

\noindent
3. 
We now check for every $u_j \in w_i|^{k_i}$ whether $u$ is a prefix of $\varphi^{-1}(u_j)|_{C(\varphi^{-1})}$, and if this is not the case, we eliminate $u_j$ from the collection of words given by $w_i|^{k_i}$. We do this for any of the $w_i \in W$ and obtain thus a collection $U$ of words $u_j$ which all have the property that $u$ is a prefix of $\varphi^{-1}(u_j)|_{C(\varphi^{-1})}$. This is precisely the finite set $U \subset F$ with the desired property
$\varphi(C_{u}^{1}) = \bigcup \limits_{u_j \in U} C_{u_{j}}^{1}$.

(The reason for this last statement is that 
the length bound, imposed  in step 2. on all $u_j \in w_i|^{k_i}$, ensures by condition (ii) above that every $C^1_{\varphi^{-1}(u_j)|_{C(\varphi^{-1})}}$ is either contained in $C^1_u$ or disjoint from the latter. Since from step 1 we know that $\varphi^{-1}(C^1_{u_j}) \subset C^1_{\varphi^{-1}(u_j)|_{C(\varphi^{-1})}}$, the same statement is true for $\varphi^{-1}(C^1_{u_j})$ replacing the $C^1_{\varphi^{-1}(u_j)|_{C(\varphi^{-1})}}$. Hence, if we eliminate in step 3 those $\varphi^{-1}(C^1_{u_j})$ from the collection which are disjoint from $C^1_u$, to determine the set $U$, then one obtains $\bigcup \limits_{u_j \in U} \varphi^{-1}(C_{u_{j}}^{1}) \subset C_{u}^{1}$ and thus  $\bigcup \limits_{u_j \in U} C_{u_{j}}^{1} \subset \varphi(C_{u}^{1})$.  

On the other hand, the inclusion $\varphi(C_{u}^{1}) \subset \bigcup \limits_{w_i \in W} C_{w_{i}}^{1} \subset \bigcup \limits_{w_i \in W} \bigcup \limits_{u_j \in w_i|^{k_i}} C_{u_{j}}^{1}$ remains true if one eliminates from the right hand term those $C_{u_{j}}^{1}$ which are disjoint from $\varphi(C_{u}^{1})$ (noting here that disjointness is preserved by the homeomorphism $\varphi$ !), which gives the converse inclusion $\varphi(C_{u}^{1}) \subset \bigcup \limits_{u_j \in U} C_{u_{j}}^{1}$.)

\end{rem}

\medskip

We'd like to point out that Llu\'is Bacardit and Ilya Kapovich have informed us that each of them observed independently the fact stated in part (1) of Remark \ref{alternative-approaches}. Furthermore, the Examples 3.9 and 3.10 in the paper 
\cite{BBC} by Berstock-Bestvina-Clay make us feel that the authors probably also had some knowledge along the lines of part (b) of Remark \ref{alternative-approaches}.
We would also like to point the reader's attention to the forthcoming paper \cite{IL}, which is in many ways a continuation of the work started here. In particular, we will treat there the question of the complexity of the algorithmic determination of the image of a given cylinder.

\medskip

We now use the results of \S 3 to define a ``dual map'' $\varphi^*$, for any automorphism $\varphi$ of $F$.  It is important, however, to always keep in mind that the definition of this map depends (heavily~!) on the choice of the basis $A$ of $F$.

\begin{defi}
\rm
Let $A$ be a basis of $F$.  For any $u \in F(A)$ we consider the finite set $U = \{ \varphi(u')|_{C(\varphi)} \mid u' \in u|^k\}$, for $k$ as in Corollary \ref{4.7}. Let $U_{\min}$ be the unique minimal set which satisfies 
$C^1_{U_{\min}} = C^1_U \,\,(= \varphi(C^1_u)$, see Proposition \ref{f12,2}).  
We define:
$$ \varphi_{A}^{*}(u)={U}_{\min}$$
Similarly, for any $U \in \mathds U$ we define $ \varphi_{A}^{*}(U)$ as the unique minimal set which defines the same cylinder as $\bigcup \limits_{u_i \in U} \varphi_{A}^{*}(u_i)$.

\end{defi}

\begin{rem}
\rm
\label{4.10}
Note that  this last definition gives directly, via Corollary \ref{4.7} and 
Proposition \ref{f12,2},
that $\varphi_{A}^{*}(u)$ does not depend on $U$ but only on $C^1_U = \varphi(C^1_u)$, and that
$\varphi(C^1_u) = \dot{\bigcup \limits_{u' \in \varphi_{A}^{*}(u)}} C^1_{u'}$.  

\end{rem}

%%%%%%%%%%%%%%%%%%%%%%

\section{Double cylinders $C^{2}_{[u,v]}$}

\begin{defi}
\rm
Let $A$ be a basis for the free group $F$.
We say that
$u,v$ are $\emph{anti-prefix}$ if  $u$ is not prefix of $v$ and $v$ is not prefix of $u$.
Similarly, we say that
$U,V \in \mathds U$ are $\emph{anti-prefix}$ if any two elements $u \in U$ and $v \in V$ are anti-prefix.

\end{defi}

\begin{rem}
\label{5.1.5}
\rm
Recall from Remark \ref{3,1} (2)
that for any $u, v \in F(A)$ the cylinders $C^1_u$ and $C^1_v$ are disjoint 
if and only if $u$ and $v$ are anti-prefix.
\end{rem}

\begin{lem}
\label{5.2}
If $u,v\in F(A)$ are anti-prefix, then $\varphi^{*}_{A}(u)$, $\varphi^{*}_{A}(v)$ are anti-prefix as well.
\end{lem}

\begin{dem}
This is a direct consequence of Remark \ref{5.1.5}, since $\varphi$ acts as homeomorphism and hence as bijection on $\partial F(A)$, so that it preserves disjointness of subsets.
\end{dem}

We now consider the Cayley graph (a tree !) $\Gamma := \Gamma(F, A)$ of the free group $F$ with respect to the basis $A$.  There is a canonical identification between the vertices of $\Gamma$ and the elements of $F$, which in turn induces a canonical identification between the boundary $\partial F$ and the set $\partial \Gamma$ of ends of $\Gamma$.  For any two $X, Y \in \partial F$ there is a well defined biinfinite reduced path $\gamma(X, Y)$ in $\Gamma$ which connects the point of $\partial \Gamma$ associated to $X$ to that associated to $Y$.

\begin{defi}
\label{5.3}
\rm
For any $u, v \in F(A)$ with $u \neq v$ we define the {\em double cylinder} $C^{2}_{[u,v]}$ as follows:
$$C^{2}_{[u,v]}= \{(X,Y)\in \partial^{2} F_{N}  \mid 
\gamma (X,Y) {\rm \, \,passes \,\, through\, \,} u {\rm \, \, and\, \,} v  {\rm \, \, (in\, \, that\, \, order)}\}$$
\end{defi}

\begin{lem}
\label{5.4}
If $u,v \in F(A)$ are anti-prefix, then one has: 
$$C^{2}_{[u,v]} = C^{1}_{u} \times C^{1}_{v}$$
\end{lem}

\begin{dem}
For $w : = u \wedge v$ (see \S 2) it follows from the assumption ``$u$ and $v$ are anti-prefix'' that $|w| < |u|$ and  $|w| < |v|$. Hence for every $(X, Y) \in C^{2}_{[u,v]}$ the geodesic $\gamma(X, Y)$ must pass (in the given order) through the points $u, w$ and $v$. In particular, it follows that $w < u < X$ and $w < v < Y$ and hence that $X \in C^{1}_{u}$ and $Y \in C^{1}_{v}$.

Conversely, for every pair $(X, Y) \in C^{1}_{u} \times C^{1}_{v}$ it follows that $w < u < X$ and $w < v < Y$, and that for $X = w \cdot X'$ and $Y = w \cdot Y'$ the biinfinite word $X'^{-1} Y'$ is reduced. Hence the geodesic $\gamma(X, Y)$ must pass (in the given order) through the points $u, w$ and $v$, which implies $(X, Y) \in C^{2}_{[u,v]}$.
\end{dem}

\begin{prop}
\label{5.5}
Let $u,v\in F(A)$ be anti-prefix. Then one has  
$$\varphi (C^{2}_{[u,v]})= \dot{\bigcup \limits_{\substack {u_{i} \in \varphi^*_{A}(u) \\ v_{j}\in \varphi^*_{A}(v)}}}  \ C^{2}_{[u_{i},v_{j}]}$$
\end{prop}

\begin{dem}
Since $u,v$ are anti-prefix, by Lemma \ref{5.4} we have $C^{2}_{[u,v]}=C^{1}_{u} \times C^{1}_{v}$,
which gives
$\varphi(C^{2}_{[u,v]})=\varphi(C^{1}_{u}) \times \varphi(C^{1}_{v})$.
By Remark \ref{4.10} we have $\varphi(C^{1}_{u}) = \dot{\bigcup \limits_{u_{i} \in \varphi^*_{A}(u)}} \ C^{1}_{u_{i}}$ and  $\varphi(C^{1}_{v}) = \dot{\bigcup \limits_{v_{j} \in \varphi^*_{A}(v)}} \ C^{1}_{v_{j}}$ and thus:
$$\varphi(C^{2}_{[u,v]})=\dot{\bigcup \limits_{u_{i} \in \varphi^*_{A}(u)}} \ C^{1}_{u_{i}} \times \dot{\bigcup \limits_{v_{j} \in \varphi^*_{A}(v)}} \ C^{1}_{v_{j}} = \dot{\bigcup \limits_{\substack{u_{i} \in \varphi^*_{A}(u) \\ {v_{j}\in \varphi^*_{A}(v)}}}} \ \big(C^{1}_{u_{i}} \times C^{1}_{v_{j}} \big)$$
By Lemma \ref{5.2} the sets $\varphi^{*}_{A}(u)$, $\varphi^{*}_{A}(v)$ are anti-prefix,
so that
by Lemma \ref{5.4} we have $C^{1}_{u_{i}} \times C^{1}_{v_{j}}  = C^{2}_{[u_{i},v_{j}]}$
for all $u_{i} \in \varphi^*_{A}(u), v_{j}\in \varphi^*_{A}(v)$,
which gives 
$$\varphi(C^{2}_{[u,v]})= \dot{\bigcup \limits_{\substack{u_{i} \in \varphi^*_{A}(u) \\ {v_{j}\in \varphi^*_{A}(v)}}}} \ C^{2}_{[u_{i},v_{j}]} \, .$$
\end{dem}

\begin{lem}
\label{5.6}
 For all $u,v,w \in F(A)$ one has  $wC^{2}_{[u,v]}= C^{2}_{[wu,wv]}$.
\end{lem}
\begin{dem}
This is a direct consequence of the definition of $C^{2}_{[u,v]}$, see Definition \ref{5.3}.
\end{dem}

Before passing to the general case of double cylinders, we need to consider the following ``small'' special cases, the proof of which follows directly from the definitions:

\begin{lem}
\label{5.7}
For any $a_i \in A$ one has:\\
$C^{2}_{[1,a_{i}]}=\dot{\bigcup \limits_{a_{j}\in A\cup A^{-1} \smallsetminus \{a_{i}\}}} \ C^{2}_{[a_{j},a_{i}]}$.\\
$C^{2}_{[1,1]}=\dot{\bigcup \limits_{a_{i}\in A\cup A^{-1}}} \ C^{2}_{[1,a_{i}]}= \dot{\bigcup \limits_{\substack {a_{j},a_{i}\in A\cup A^{-1} \\ a_{i} \neq a_{j}}}} \ C^{2}_{[a_{j},a_{i}]}$.
\qed
\end{lem}

\begin{prop}
\label{5.8}
For any two distinct $u,v\in F(A)$ there exist finite computable sets $U,V \subset F(A)$ such that
$$\varphi (C^{2}_{[u,v]})=\dot{\bigcup \limits_{\substack {u_{i} \in U \\ v_{j}\in V}}}  \ C^{2}_{[u_{i},v_{j}]}$$
\end{prop}

\begin{dem}
If $u$ and $v$ are anti-prefix, then Proposition \ref{5.5} gives the desired statement (and furthermore a precise description of the sets $U$ and $V$).

Otherwise, one has $u \leq v$ or $v \leq u$, and if $\big{|}|u|-|v|\big{|} \geq 2$ we can find some $w \in F(A)$ with $u < w < v$ or $v < w < u$. Hence  Lemma \ref{5.6} allows us to replace $u$ by $w^{-1} u$ and $v$ by $w^{-1} v$, which reduces this case to the one treated in the previous paragraph.

Finally, if $\big{|}|u|-|v|\big{|} \leq 1$ we can first again apply Lemma \ref{5.6} to achieve that $u = 1$ or $v = 1$.  But then Lemma \ref{5.7} brings us again back to the case treated in the first paragraph.
\end{dem}

\begin{rem}
\label{pas-encore-ecrit}
\rm
From the arguments given in the last proof one can derive the following improvement of Proposition \ref{5.5}:

For any two distinct $u,v\in F(A)$ 
(i.e. without supposing that they are anti-prefix) 
one has:
$$\varphi (C^{2}_{[u,v]})=\dot{\bigcup \limits_{\substack {u_{i} \in \varphi(v) \varphi^*_{A}(v^{-1}u) \\ v_{j}\in \varphi(u)\varphi^*_{A}(u^{-1} v)}}}  \ C^{2}_{[u_{i},v_{j}]}$$
\end{rem}


\begin{thebibliography}{ABHS05}


\bibitem{BBC} 
J.~Behrstock, 
M.~Bestvina, and M.~Clay,
{\em Growth of intersection numbers for free group automorphisms.} arXiv:0806.4975 


\bibitem{BH}
M.~Bridson, A.~Haefliger, {\em Metric spaces of nonÐpositive curvature.} Springer, New York
1999

\bibitem{MR916179}
D.~Cooper,
{\em Automorphisms of free groups have finitely
generated fixed point
  sets.}
 J. Algebra {\bf 111}, 453--456 (1987)

\bibitem{FI}
F.~Ibrahim, {\em Cylindres, multi-cylindres et leur images sous l'action de $\Aut(F_n)$.}
Ph.D.-thesis, Universit\'e d'Aix-Marseille, 2012

\bibitem{IL}
F.~Ibrahim, M.~Lustig, 
{\em Dual automorphisms of free groups.} Preliminary preprint, 2012

\end{thebibliography}
\end{document}